\documentclass[11pt]{amsart}
\usepackage{amsmath,amssymb,amsthm,amsfonts}
\usepackage{mathrsfs,eufrak,euscript}
\usepackage{enumerate}
\theoremstyle{plain}
\newtheorem{thm}[subsection]{Theorem}
\newtheorem{prop}[subsection]{Proposition}
\newtheorem{lemma}[subsection]{Lemma}
\newtheorem{eg}[subsection]{Example}
\newtheorem{eg's}[subsection]{Examples}

\newtheorem{defn}[subsection]{Definition}
\newtheorem{rmk}[subsection]{Remark}

\newtheorem{note}[subsection]{Note}

\title[Continuous functional calculus]{Continuous functional calculus for  quaternionic bounded normal operators}
\author[G. Ramesh {\protect \and} P. Santhosh Kumar]{G. Ramesh {\protect
		\and} P. Santhosh Kumar}
\address{G. Ramesh, Department of Mathematics, Kandi, Sangareddy, Telangana, India 502285.}
\email{rameshg@iith.ac.in}

\address{P. Santhosh Kumar, Department of Mathematics, Kandi, Sangareddy, Telangana, India 502285.}
\email{ma12p1004@iith.ac.in}

\begin{document}
	\maketitle		
	\begin{abstract}
		In this article we give an approach to define continuous functional calculus for bounded quaternionic normal operators defined on a right quaternionic Hilbert space.
		
	\end{abstract}
	\section{Introduction and Preliminaries}
	Several authors discussed the significance of functional calculus, in quaternionic Hilbert spaces \cite{Adler,Birkhoff,Emch, Emch1, Finkelstein1, Horwitz}. A special class of functions with quaternion domain named {\it slice regular functions} is introduced by Gentili and Struppa \cite{Gentili1, Gentili2}. This class of functions are the appropriate generalization of standard holomorphic functions. The theory of slice hyperholomorphic functions, related S-functional calculus, Riesz-Dunford functional calculus for bounded operators are studied in  \cite{Ghiloni}, whereas the S-functional calculus for closed densely defined operators can be found in \cite{Gantner} and the $H^\infty$- functional calculus for $n$- tuple of non - commuting quaternionic operators in  \cite{Alpay3}.
	
	There are several versions of the spectral theorem in quaternionic setting in the literature. The spectral theorem that deals with the integral representation of a quaternionic normal operator is given by Viswanath  \cite{viswanath}.  The author proved the existence of  spectral measure through the symplectic image and as a consequence, obtained the  Cartesian decomposition of a normal operator in a quaternionic Hilbert space. Sushama Agarwal and S. H. Kulkarni    \cite{Sushama} proved the spectral theorem for normal operators on  real Hilbert spaces by exploiting real Banach algebra techniques, and deduced the quaternionic version from this. The spectral theorem for quaternionic unitary operators is proved in \cite{Sharma}. The same result was obtained by  Alpay et al. using the  notion of spherical spectrum and the quaternionic version of the Herglotz theorem in \cite{Alpay1}, later generalized it to the case of unbounded normal operators  \cite{Alpay2}. A similar result using a different approach can also be found in \cite{Ghiloni2}.
	
	In this article we give an approach to define continuous functional calculus for bounded quaternionic normal operators, and as a consequence we deduce the integral representation (spectral theorem) of such operators.  First we define the functional calculus for a subclass of continuous functions by using the classical functional calculus from \cite{Zhu}, then we extend the same to quaternion valued continuous functions.	

	We organize this article in three sections. In the first section, we give introduction to the ring of quaternions, quaternionic Hilbert spaces and recall some of the basic results which we need for our purpose.
	
	In the second section, we give an approach to continuous functional calculus for bounded quaternionic normal operators. In the final section, we establish the integral representation of continuous functions of bounded quaternionic normal operators.
	\subsection*{Quaternions:}
	
	Let $1, i, j$ and $k$ be the vectors of the canonical orthonormal basis for the Euclidean space $\mathbb{R}^{4}$, which satisfies the following property:
	\begin{equation*}
	i^{2}= j^{2}=k^{2} = -1 = i\cdot j \cdot k.
	\end{equation*}
	Let $\mathbb{H}$ denote the 4-dimensional real algebra consisting of elements, called real quaternions, of the form
	\begin{equation}\label{eq:expressionq}
	q = 	q_{0}+q_{1}i+q_{2}j+q_{3}k,
	\end{equation}
	where $q_{\ell}\in \mathbb{R}$ for $\ell = 0,1,2,3.$
	It is easy to see that $\mathbb{H}$ is a division ring (skew field) of all real quaternions.
	The real and imaginary parts of $q$ in Equation (\ref{eq:expressionq}) are given by,$
	\text{re}(q) = q_{0} \; ; \; \text{im}(q)= q_{1}i+q_{2}j+q_{3}k$. The conjugate of $q$, denoted by $\overline{q}$, is defined by $ \overline{q} = q_{0} - q_{1} i - q_{2} j - q_{3} k$.
	The modulus of $q$ is given by $|q|^{2} = \sum\limits_{\ell = 0}^{3} q_{\ell}^{2}$. This defines the norm on $\mathbb{H}$. If $ q \in \mathbb{H}\setminus \{0\}$, then $q$ is invertible and the inverse is given by $\displaystyle q^{-1} = \frac{\overline{q}}{ |q|^{2}}$.
	The imaginary unit sphere in $\mathbb{H}$, denoted by $\mathbb{S}$, is defined by $	\mathbb{S}:= \big\{ q \in \mathbb{H}:\; \overline{q}=-q,\; |q| = 1 \big\}$.
	For each $m \in \mathbb{S}$, the real subalgebra  $\mathbb{C}_{m}:= \{\alpha + m  \beta \; ;\; \alpha, \beta \in \mathbb{R} \}$ is called the slice of $\mathbb{H}$. In fact, $\mathbb{C}_{m}$ is a field, for all $m \in \mathbb{S}$. We denote the upper half plane of $\mathbb{C}_{m}$ by $\mathbb{C}_{m}^{+}:= \big\{\alpha + m \beta \; |\; \alpha \in \mathbb{R}, \beta \geq 0\big\}$.
	
	Note that for each $m \in \mathbb{S}$, the slice $\mathbb{C}_{m}$ is isomorphic (algebra) to $\mathbb{C}$ through the map $ \alpha + m\beta \mapsto \alpha+i\beta $. So, all the results in the theory of complex Hilbert spaces holds true for $\mathbb{C}_{m}$- Hilbert spaces.	
	
	Here we list out some of the properties of quaternions (see \cite{Ghiloni} for details), which we need later. Let $p,q \in \mathbb{H}$. Then
	
	\begin{enumerate}
		\item $ |p.q| = |p|.|q|$ and $|\overline{p}| = |{p}|$.
		\item  $ q \in \mathbb{C}_{m} $ for some $m \in \mathbb{S}$ if and only if $ q.\lambda = \lambda.q,$ for every $\lambda \in \mathbb{C}_{m}$.
		\item $\mathbb{C}_{m} \cap \mathbb{C}_{n} = \mathbb{R}$, for $ m \neq \pm n \in \mathbb{S}$. Moreover, $\mathbb{H} = \bigcup\limits_{m \in \mathbb{S}} \mathbb{C}_{m}$.
		

	\end{enumerate}
	We recall that there is an equivalence relation on $\mathbb{H}$ defined as, $p\sim q$ if and only if $ p = s^{-1}qs $, for some $s \in \mathbb{H}\setminus \{0\}$.
	The equivalence class of $ q$ is given by
	\begin{equation*}
	[q]: = \big\{ p \in \mathbb{H}: p \sim q \big\} = \big\{ p \in \mathbb{H}: \text{re}(p) = \text{re} (q)\; \text{and}\; |\text{im}(p)| = |\text{im}(q)| \big\}.
	\end{equation*}
	\begin{defn}\label{circularization} \cite[Equation 5.19]{Ghiloni}
		Let $m \in \mathbb{S}$ and $\mathcal{K}$ be a subset of $\mathbb{C}_{m}$. Then the circularization of $\mathcal{K}$, denoted by $\displaystyle\Omega_{\mathcal{K}}$, is defined  by
		\begin{equation*}
		\Omega_{\mathcal{K}}:= \big\{\alpha + m^{\prime} \beta \colon \ \alpha, \beta \in \mathbb{R}, \alpha + m \beta \in \mathcal{K}, m^{\prime} \in \mathbb{S}  \big\}.
		\end{equation*}
	\end{defn}
	
	\begin{defn}
		Let $ \mathcal{H} $ be a right $ \mathbb{H}- $ module. A map $ \left\langle \cdot  | \cdot \right\rangle \colon \mathcal{H} \times \mathcal{H} \to \mathbb{H}$ is called an inner product on $\mathcal{H}$ if it satisfies the following three properties:
		\begin{enumerate}
			\item  $\left\langle x | x\right\rangle \geq 0$, for all $x \in \mathcal{H}$ and $\left\langle x | x\right\rangle = 0$ if and only if $x = 0$.
			\item $\left\langle x | yp+zq\right\rangle = \left\langle x| y \right\rangle p + \left\langle x | z\right\rangle q, $ for all $x,y,z \in \mathcal{H}$ and $p,q \in \mathbb{H} $.
			\item $\left\langle x | y\right\rangle = \overline{\left\langle y | x\right\rangle}$ , for all $x,y \in \mathcal{H}$.
		\end{enumerate}
		Define $ \| x \| := {\left\langle x| x \right\rangle}^{\frac{1}{2}}, $ for every $ x \in \mathcal{H}.$  Then $\|\cdot\|$ is a norm on $\mathcal{H}$. If the normed space $ (\mathcal{H}, \|\cdot\|) $ is complete, then  $\mathcal{H}$ is called a right quaternionic Hilbert space.
	\end{defn}
	\begin{note}
		Throughout this article $\mathcal{H}$ denotes a right quaternionic Hilbert space.	
		If $x,y \in \mathcal{H}$, then the following polarization identity (see \cite[Proposition 2.2]{Ghiloni} for details) holds:
		\begin{equation}\label{polarization}
		4 \left\langle x | y \right\rangle= \sum_{l=1,i,j,k}\Big(\|xl+y\|^{2}-\|xl-y\|^{2}\Big)\cdot l.
		\end{equation}
	\end{note}
\begin{eg}
%
 Let $m \in \mathbb{S}$ and $\Omega \subseteq \mathbb{C}_{m}$. Let $\mu$ be a positive $\sigma$ - additive measure on $\Omega$, then
		\begin{equation*}
		L^{2}(\Omega; \mathbb{H}; \mu) : = \left\{f \colon \Omega \to \mathbb{H} \; | \int\limits_{\Omega}|f(x)|^{2}d\mu(x) < \infty\right\}
		\end{equation*}
		is a right quaternionic Hilbert space with the scalar multiplication and inner product defined, respectively, by
		\begin{align*}
		(f\cdot p)(x) &= f(x)\cdot p, \; \text{for all}\;  f \in L^{2}(\Omega; \mathbb{H}; \mu),\; p \in \mathbb{H}\; ;\;\\
		\left\langle f \;\big|\; g\right\rangle &= \int\limits_{\Omega} \overline{f(x)}\cdot g(x) \ d\mu(x), \; \text{for all}\; f,g \in L^{2}(\Omega; \mathbb{H}; \mu).
		\end{align*}
\end{eg}

	\begin{defn}
		Let $\mathcal{S} \subseteq \mathcal{H}$. The orthogonal complement of $\mathcal{S}$, denoted by $\mathcal{S}^{\bot}$, is defined by $
		{\mathcal{S}}^{\bot} : = \big\{ x \in \mathcal{H} \ | \ \left\langle x | y \right\rangle = 0, \; \text{for all} \; y \in {\mathcal{S}} \big\}$.
	\end{defn}
	\begin{thm}
		Let $\mathcal{N}$ be a subset of $\mathcal{H}$ such that , for $z, z^\prime \in \mathcal{N}$, $\left\langle z | z^\prime\right\rangle = 0$ if $z \neq z^\prime$ and $\left\langle z | z\right\rangle = 1$. Then  the set $\mathcal{N}_{x}:= \big\{z \in \mathcal{N}:\; \left\langle z|x\right\rangle \neq 0\big\}$	is countable, for all $x \in \mathcal{H}$.
	\end{thm}
Now we define Hilbert basis of a right quaternionic Hilbert space (see \cite[Proposition 2.5]{Ghiloni} for details).
\begin{defn}
	A subset $\mathcal{N}$ of  $\mathcal{H}$ is said to be a Hilbert basis of $\mathcal{H}$ if, for every $z,z^{\prime} \in \mathcal{N}$, we have $\left\langle z | z^{\prime}\right\rangle = \delta_{z,z^{\prime}}$ and $\left\langle x | y \right\rangle = \sum\limits_{z \in \mathcal{N}} \left\langle x | z\right\rangle \left\langle z | y \right\rangle$ for all $x,y \in \mathcal{H}$.
\end{defn}
\begin{rmk}
	Every quaternionic Hilbert space $\mathcal{H}$ admits a Hilbert basis $\mathcal{N}$ (see \cite[Proposition 2.6]{Ghiloni}), and every $x \in \mathcal{H}$ can be uniquely decomposed as follows:
	\begin{equation*}
	x = \sum\limits_{z \in \mathcal{N}} z \left\langle z | x \right\rangle.
	\end{equation*}
\end{rmk}
	\begin{defn}\cite[Definition 2.9]{Ghiloni}
		A map  $T \colon \mathcal{H} \to \mathcal{H} $  is said to be right $\mathbb{H}$- linear or quaternionic operator, if
		$T(x+y) = Tx + Ty $ and $T(x \cdot q)= Tx \cdot q$, for all $ x,y \in \mathcal{H}, q \in \mathbb{H}$. Moreover, if there exist $M>0$ such that $\|Tx\| \leq M \|x\|$, for all $x \in \mathcal{H}$, then $T$ is called bounded or continuous. In this case, the norm of $T$, defined by
		\begin{equation*}
		\|T\| = \sup \Big\{\|Tx\|:\; x \in \mathcal{H}, \|x\|=1\Big\},
		\end{equation*}
		is finite.
	\end{defn}
We denote the class of all bounded quaternionic operators on $\mathcal{H}$ by $\mathcal{B}(\mathcal{H})$.
	
	
	\begin{defn}\cite[Definition 2.2]{Alpay2}
		Let $ T \in \mathcal{B}(\mathcal{H})$. Then there exists unique $T^{\ast} \in \mathcal{B}(\mathcal{H})$ such that
		\begin{equation*}
		\left\langle x | Ty\right\rangle = \left\langle T^{*}x | y\right\rangle, \; \text{for all}\; x,y \in \mathcal{H}.
		\end{equation*}
		This operator $ T^{*} $ is called the adjoint of $T.$
	\end{defn}
	%
	\begin{defn}\cite{Ghiloni}
		Let $T \in \mathcal{B}(\mathcal{H})$. Then $T$ is said to be normal, if $ T^{*}T = TT^{*}$,	self-adjoint if $T^{*}=T$. We say $T$ to be positive if $T=T^{*}$ and $ \left\langle x | Tx \right\rangle \geq 0,$ for all $ x\in \mathcal{H}$, anti self-adjoint if $T^{*}=-T$ and unitary if $ TT^{*}=T^{*}T= I$.
		\end{defn}
	\begin{defn} \cite[Section 3.1]{Ghiloni} \label{Leftmultiplication}
		Let $\mathcal{N}$ be a Hilbert basis of $\mathcal{H}$. Then the left multiplication, induced by $\mathcal{N}$, is defined by a map $(q,x) \in \mathbb{H} \times \mathcal{H} \mapsto q \cdot x \in \mathcal{H}$ that is
		\begin{equation}\label{leftmulti}
		q \cdot x = \sum\limits_{z \in \mathcal{N}} z \cdot q \left\langle z | x \right\rangle, \; \text{for every} \; q \in \mathbb{H} \; \text{and} \; x \in \mathcal{H}.
		\end{equation}
	\end{defn}
	\begin{note} Let $q \in \mathbb{H}$. Then the map $L_{q} \colon \mathcal{H} \to \mathcal{H}$ defined by
		\begin{equation*}
		L_{q}(x):= q \cdot x, \; \text{for all} \; x \in \mathcal{H}
		\end{equation*}
		is bounded by Definition \ref{Leftmultiplication}. Moreover, $L_{q}^{*} = L_{\overline{q}}$ and $\|L_{q}\| = |q|$. We see that  $L_{q}$ is anti self-adjoint and unitary if and only if $q \in \mathbb{S}$. 
	\end{note}
	\begin{lemma}\cite[Lemma 4.1]{Ghiloni}\label{directsum}
		Let $\left\langle \cdot | \cdot \right\rangle \colon \mathcal{H} \times \mathcal{H} \to \mathbb{H}$ be an inner product on $\mathcal{H},\;m \in \mathbb{S}$ and let $J$ be an anti self-adjoint unitary operator on $\mathcal{H}$. Define $ \mathcal{H}^{Jm}_{\pm}= \big\{x \in \mathcal{H}: J(x)= \pm x\cdot m\big\}$.
		Then
		\begin{enumerate}
			\item $\mathcal{H}^{Jm}_{\pm} \neq \{0\}$ and the restriction of the inner product $\left\langle \cdot | \cdot \right\rangle$ to $\mathcal{H}^{Jm}_{\pm}$ is $\mathbb{C}_{m}$- valued. Therefore $\mathcal{H}^{Jm}_{\pm}$ is $\mathbb{C}_{m}-$ Hilbert space, called the slice Hilbert space of $\mathcal{H}$.
			\item $\mathcal{H} = \mathcal{H}^{Jm}_{+} \oplus \mathcal{H}^{Jm}_{-}$.
		\end{enumerate}
	\end{lemma}
	We denote the class of all bounded $\mathbb{C}_{m}$- linear operators on $\mathcal{H}^{Jm}_{+}$ by $\mathcal{B}(\mathcal{H}^{Jm}_{+})$.
	\begin{rmk}\cite[Proposition 3.8(f)]{Ghiloni} \label{J}
		Let $m \in \mathbb{S}$. If $\mathcal{N}$ is a Hilbert basis of $\mathcal{H}^{Jm}_{+}$, then $\mathcal{N}$ is also a Hilbert basis of $\mathcal{H}$ and
		\begin{equation*}
		J(x) = \sum\limits_{z\in \mathcal{N}} z \cdot m \left\langle  z | x \right\rangle = L_{m}(x), \; \text{for all} \; x \in \mathcal{H}.
		\end{equation*}
		This implies that $J = L_{m}$.
	\end{rmk}
We recall the notion of the spherical spectrum of quaternionic operators.
	\subsection*{Spherical spectrum} \cite[Definition 4.1]{Ghiloni}
	Let  $ T \in \mathcal{B}(\mathcal{H})$  and $ q \in \mathbb{H}.$ Define  $\Delta_{q}(T) \colon \mathcal{H} \to \mathcal{H}$ by
	\begin{equation*}
	\Delta_{q}(T):= T^{2}-T(q+\overline{q})+I\cdot|q|^{2}.
	\end{equation*}
	The spherical spectrum of $T$, denoted by $\sigma_{S}(T)$, is defined by
	\begin{equation*}
	\sigma_{S}(T):= \Big\{q \in \mathbb{H}:\; \Delta_{q}(T)\;\text{is not invertible in}\; \mathcal{B}(\mathcal{H})\Big\}.
	\end{equation*}
	Note that $\sigma_{S}(T)$ is a non-empty compact subset of $\mathbb{H}$.
	
	It is proved that every densely defined linear operator on a slice Hilbert space can be extended uniquely to a densely defined right $\mathbb{H}$- linear operator on a quaternionic Hilbert space, and the converse is true under certain conditions \cite[Proposition 3.11]{Ghiloni}. The same result is true for bounded operators. Here we quote this result for bounded operators.
	\begin{prop} \label{extension1}
		For every \ $\mathbb{C}_{m}$- linear operator $T \colon  \mathcal{H}^{Jm}_{+} \to \mathcal{H}^{Jm}_{+}$,  there exists a unique quaternionic operator  $ \widetilde{T} \colon\mathcal{H} \to \mathcal{H}$
		such that $ \widetilde{T}(x) = T(x), $ for every $x \in \mathcal{H}^{Jm}_{+}.$ Moreover, the following facts holds:
		\begin{enumerate}
			\item If $T \in \mathcal{B}(\mathcal{H}^{Jm}_{+})$, then $\widetilde{T}\in \mathcal{B}(\mathcal{H})$ and $\|\widetilde{T}\| = \|T\|$.
			\item $J\widetilde{T} = \widetilde{T} J$.
		\end{enumerate}
		On the other hand, let $V \colon \mathcal{H} \to \mathcal{H}$ be a bounded quaternionic operator. Then $ V = \widetilde{U} $, for a unique bounded $\mathbb{C}_{m}-$ linear operator $ U \colon \mathcal{H}^{Jm}_{+} \to \mathcal{H}^{Jm}_{+} $ if and only if  $JV = VJ$.	Furthermore,
		\begin{enumerate}
			\item  $\big(\widetilde{T}\big)^{*} = \widetilde{T^{*}}$.
			\item \label{extnmulti}If $S \colon  \mathcal{H}^{Jm}_{+} \to \mathcal{H}^{Jm}_{+}$ is $\mathbb{C}_{m}$- linear, then $\widetilde{ST} = \widetilde{S} \widetilde{T}$.
			\item \label{extninverse}If $S$ is the inverse of $T$, then $\widetilde{S}$ is the inverse of $\widetilde{T}$.
		\end{enumerate}
	\end{prop}
\begin{rmk}\label{rmkforboundedJ}
	In particular, if $T \in \mathcal{B}(\mathcal{H})$ is normal operator, then by  \cite[Theorem 5.9]{Ghiloni}, there exist an anti self-adjoint and unitary $J \in \mathcal{B}(\mathcal{H})$ such that $TJ = JT $.
\end{rmk}

	\begin{note}\label{note:standardspectrum}
		If $\widetilde{T}$ is the extension of a bounded $\mathbb{C}_{m}$- linear operator $T$ on $\mathcal{H}^{Jm}_{+}$, then by \cite[Corllary 5.13]{Ghiloni}, the spectrum of $T$ is give by
		\begin{equation*}
		\sigma(T) = \sigma_{S}(\widetilde{T})\cap \mathbb{C}_{m}^{+}.
		\end{equation*}
	\end{note}

	\section{Continuous functional calculus}
	
	In this section we give an approach to define continuous functional calculus for bounded quaternionic normal operators. We recall that the circularization of a non-empty set $\mathcal{K}$ of $\mathbb{C}_{m} (m \in \mathbb{S})$ is given by
	\begin{equation*}
	\Omega_{\mathcal{K}} = \big\{\alpha + m^{\prime} \beta : \alpha, \beta \in \mathbb{R}, \alpha + m\beta \in \mathcal{K}, m^{\prime}  \in \mathbb{S}\big\}.
	\end{equation*}
	Let $C(\Omega_{\mathcal{K}}, \mathbb{H})$ denote the class of all $\mathbb{H}$- valued continuous functions on ${\Omega}_{\mathcal{K}}$ i.e.,
	\begin{equation*}
	C(\Omega_{\mathcal{K}}, \mathbb{H}) = \big\{f\colon \Omega_{\mathcal{K}}\to \mathbb{H}: f \; \text{is continuous}\big\}.
	\end{equation*}
	It is a real algebra with the addition and the multiplication defined respectively by, $
	(f+g)(q) = f(q)+g(q) \; \text{and}\; (f\cdot g)(q) =f(q)g(q)$
	for all $f,g \in {C}(\Omega_{\mathcal{K}}, \mathbb{H})$ and $q \in \Omega_{\mathcal{K}}$. The scalar multiplication is defined by
	\begin{equation*}
	(rf)(q) = r f(q) \; \text{and}\; (fr)(q) = f(q)r, \; \text{for all} \; f \in {C}(\Omega_{\mathcal{K}}, \mathbb{H}), r \in \mathbb{R}.
	\end{equation*}
	the conjugate of $f$, denoted by $\overline{f}$, is defined by
	\begin{equation*}
	\overline{f}(q) = \overline{f(q)}, \; \text{for all} \; q \in \mathbb{H}.
	\end{equation*}
	 If  $f\in {C}(\Omega_{\mathcal{K}}, \mathbb{H})$ is bounded, then
	\begin{equation*}
	\|f\|_{\infty} = \sup \big\{|f(q)| : \; q \in \Omega_{\mathcal{K}}\big\}.
	\end{equation*}
	For $m \in \mathbb{S}$, we introduce a subclass of $C(\Omega_{\mathcal{K}}, \mathbb{H})$,  that is
	\begin{equation*}
	{C}_{\mathbb{C}_{m}}(\Omega_{\mathcal{K}}, \mathbb{H}): = \big\{f \in {C}(\Omega_{\mathcal{K}}, \mathbb{H})  : f(\mathcal{K}) \subseteq \mathbb{C}_{m}\big\}.
	\end{equation*}
	It is immediate to see that ${C}_{\mathbb{C}_{m}}(\Omega_{\mathcal{K}}, \mathbb{H})$ is a real subalgebra of ${C}(\Omega_{\mathcal{K}}, \mathbb{H})$.
	\begin{lemma}\label{invconti}
		Let $m \in \mathbb{S}$ and $\mathcal{K} \subseteq \mathbb{C}_{m}$. If  $f \in {C}_{\mathbb{C}_{m}}(\Omega_{\mathcal{K}}, \mathbb{H})$, then the map $f_{+} \colon \mathcal{K} \to \mathbb{C}_{m}$ defined by
		\begin{equation*}
		f_{+}(z) = f(z) , \; \text{for all}\; z \in \mathcal{K}
		\end{equation*}
		is continuous.
	\end{lemma}
	\begin{proof} We show that  $f_{+}$ is well defined. Let $z_{1},z_{2} \in \mathcal{K}$ be such that $z_{1} = z_{2}$. Then by  the definition of $f_{+}$, we have
		\begin{align*}
		f_{+}(z_{1}) = f(z_{1}) = f(z_{2}) = f_{+}(z_{2}).
		\end{align*}
		Now we show that $f_{+}$ is continuous. Let $(z_{\ell}) \subset \mathcal{K}$ be a sequence such that $(z_{\ell})$ converges to $z_{0}$ in $\mathcal{K}$, as $\ell\to \infty$. Since $\mathcal{K} \subset \Omega_{\mathcal{K}}$ and $f$ is continuous, we see that
		\begin{equation*}
		f_{+}(z_{\ell}) = f(z_{\ell}) \longrightarrow f(z_{0}) = f_{+}(z_{0}),
		\end{equation*}
		as $\ell \to \infty$.	This shows that  $f_{+}$ is continuous.
	\end{proof}
	
	By Lemma \ref{invconti}, it is clear that if  $f \in {C}_{\mathbb{C}_{m}}(\Omega_{\mathcal{K}}, \mathbb{C}_{m})$, then $f_{+} \in C({\mathcal{K}}, \mathbb{C}_{m})$, the class of all  $\mathbb{C}_{m}$- valued continuous functions on $\mathcal{K}$.

	\begin{note} \label{noteforT+}
		Let $m,n \in \mathbb{S}$ with $mn = - nm$ and $T \in \mathcal{B}(\mathcal{H})$ be normal. By Remark \ref{rmkforboundedJ}, there exists an anti self-adjoint and unitary operator $J \in \mathcal{B}(\mathcal{H})$ such that $JT = TJ$. Then  by Theorem \ref{extension1}, there exist a unique $\mathbb{C}_{m}$- linear operator  $T_{+} \in \mathcal{B}( \mathcal{H}^{Jm}_{+})$  such that $T = \widetilde{T_{+}}$.
		%
	\end{note}
	\begin{rmk}\label{rmkforspT+}
		Let $T \in \mathcal{B}(\mathcal{H})$ be normal and $T_{+}$ be as in Note \ref{noteforT+}. Then
		by Note \ref{note:standardspectrum}, we have
		\begin{equation*}
		\sigma_{S}(T) = \Omega_{\sigma(T_{+})}.
		\end{equation*}	
	\end{rmk}
	We refer \cite{Zhu} for the functional calculus in complex Hilbert spaces.
	\begin{defn}\label{Definitioncmf}
		Let $m,n \in \mathbb{S}$ with $mn = - nm$. Let $T \in \mathcal{B}(\mathcal{H})$ be normal and $T_{+} \in \mathcal{B}(\mathcal{H}^{Jm}_{+})$ be such that $\widetilde{T_{+}} = T$. If $f \in {C}_{\mathbb{C}_{m}}(\sigma_{S}(T), \mathbb{H})$, then by Lemma \ref{invconti}, $f_{+}: = f|_{\sigma(T_{+})}$ is continuous. By the continuous functional calculus for normal operators on a complex Hilbert space, $f_{+}(T_{+})$ is well defined and it is a normal operator. We define
		\begin{equation} \label{defofF_calculus}
		f(T):= \widetilde{f_{+}(T_{+})}.
		\end{equation}	
	\end{defn}
	Since $f_{+} \in C(\sigma(T_{+}), \mathbb{C}_{m})$ and  $T_{+} \in \mathcal{B}(\mathcal{H}^{Jm}_{+})$, by continuous functional calculus in complex Hilbert spaces \cite{Zhu}, the operator $f_{+}(T_{+})$ is  well-defined. Thus from Equation (\ref{defofF_calculus}) and by Proposition \ref{extension1}, $ f(T)$ is well-defined.
	\begin{lemma}
		Let $T \in \mathcal{B}(\mathcal{H})$ be normal and $m \in \mathbb{S}$. If $f \in {C}_{\mathbb{C}_{m}}(\sigma_{S}(T), \mathbb{H})$, then $f(T)^{*} = \overline{f}(T)$.
	\end{lemma}
	\begin{proof}
		By Definition \ref{Definitioncmf}, we have
		\begin{align*}
		f(T)^{*} & = \widetilde{f_{+}(T_{+})^{*}} = \widetilde{\overline{f_{+}}(T_{+})} = \overline{f}(T),
		\end{align*}
		where $f_{+}=f|_{\sigma(T_{+})}$. This completes the proof.
	\end{proof}
	\begin{rmk}\label{boundedapprx}
		Let $T \in \mathcal{B}(\mathcal{H})$ be normal and $m \in \mathbb{S}$. If $f \in {C}_{\mathbb{C}_{m}}(\sigma_{S}(T), \mathbb{H})$ and $f_{+}:=f|_{\sigma(T_{+})}$, then by Weierstrass approximation theorem there exist a sequence of polynomials $(p_{\ell})$ with coefficients from $\mathbb{C}_{m}$ such that $\|p_{\ell}- f_{+}\|_{\infty} \longrightarrow 0$, as $\ell \to \infty$. In fact, $f_{+}(T_{+})$ is defined by
		\begin{equation*}
		f_{+}(T_{+}):= \lim\limits_{\ell \to  \infty} p_{\ell}(T_{+})
		\end{equation*}
		This implies that
		\begin{equation*}
		f(T): = \lim\limits_{\ell \to \infty}p_{\ell}(\widetilde{T_{+}}) = \lim \limits_{\ell \to \infty}p_{\ell}(T).
		\end{equation*}
	\end{rmk}
	Now we generalize functional calculus to the class of all $\mathbb{H}$- valued continuous functions on $\sigma_{S}(T)$.  Let $m,n \in \mathbb{S}$ with $mn = - nm$. Then every $q \in \mathbb{H}$ is uniquely expressed as follows:
	\begin{equation}\label{decompositionofq}
	q = \frac{q+ \overline{q}}{2} - \frac{(qm+ \overline{qm})m}{2} -\frac{(qn+ \overline{qn})n}{2}-\frac{(qmn+ \overline{qmn})mn}{2}.
	\end{equation}
	In the view of Equation (\ref{decompositionofq}), every $f \in {C}(\sigma_{S}(T), \mathbb{H})$ can be uniquely decomposed as the following:
	\begin{equation}\label{decompositionoff}
	f(q) = f_{0}(q) + f_{1}(q)m + f_{2}(q)n + f_{3}(q)mn, \; \text{for all} \; q \in \sigma_{S}(T),
	\end{equation}
	where
	\begin{align*}
	f_{0}(q) = \frac{f(q)+\overline{f(q)}}{2}, \; & \; f_{1}(q) = -\frac{f(q)m+  \overline{f(q)m}}{2}, \\
	f_{2}(q) = - \frac{f(q) + \overline{f(q)n}}{2}\; &\; \text{and}\; f_{3}(q) = -\frac{f(q)+  \overline{f(q)mn}}{2}.
	\end{align*}
	Since $f$ is continuous, we see  that $f_{\ell} \in C(\sigma_{S}(T), \mathbb{R})$, for $\ell = 0,1,2,3$.
	
	
	\begin{lemma}\label{decomposition}
		Let $T \in \mathcal{B}(\mathcal{H})$ be normal and $m,n \in \mathbb{S}$ with $mn = - nm$. If $f \in C(\sigma_{S}(T), \mathbb{H})$, then there exist unique $F_{1}, F_{2} \in C(\sigma_{S}(T), \mathbb{C}_{m})$ such that
		\begin{equation}\label{decompositioneq}
		f(q) = F_{1}(q) + F_{2}(q) \cdot n, \; \text{for all}\; q \in \sigma_{S}(T).
		\end{equation}
	\end{lemma}
	\begin{proof}
		Since $f \in C(\sigma_{S}(T), \mathbb{H})$, then by Equation (\ref{decompositionoff}), we have
		\begin{equation*}
		f(q) = f_{0}(q) + f_{1}(q)m+f_{2}(q)n+f_{3}mn, \; \text{for all}\; q \in \sigma_{S}(T).
		\end{equation*}
		Define $F_{1}, F_{2}$ on $\sigma_{S}(T)$ by
		\begin{equation*}
		F_{1}(q) = f_{0}(q) + f_{1}(q)m \; ; \;F_{2}(q) = f_{2}(q) + f_{3}(q)m, \; \text{for all} \; q \in \sigma_{S}(T).
		\end{equation*}
		Since $f_{\ell} \in C(\sigma_{S}(T), \mathbb{R})$, we have $F_{1}, F_{2} \in C(\sigma_{S}(T), \mathbb{C}_{m})$. Moreover,
		\begin{equation*}
		f(q) = F_{1}(q) + F_{2}(q) \cdot n, \; \text{for all}\; q \in \sigma_{S}(T).
		\end{equation*}
		Let $G_{1}, G_{2} \in C(\sigma_{S}(T), \mathbb{C}_{m})$ be such that $f(q) = G_{1}(q)+ G_{2}(q) \cdot n$, for all $q \in \sigma_{S}(T)$. Then
		\begin{equation*}
		[F_{1}(q)-G_{1}(q)] + [F_{2}(q) - G_{2}(q)] \cdot n = 0, \; \text{for all}\; q \in \sigma_{S}(T).
		\end{equation*}
		This implies that $F_{1}(q) = G_{1}(q)$ and $F_{2}(q) = G_{2}(q)$. Hence the decomposition of $f$ as in Equation (\ref{decompositioneq}) is unique.
	\end{proof}
	\begin{note} \label{noteforJ}
		Let $m,n \in \mathbb{S}$ with $mn = -nm$, and $J \in \mathcal{B}(\mathcal{H})$ be anti self-adjoint and unitary. Then by \cite[Proposition 3.8]{Ghiloni}, there exists a scalar multiplication $q \mapsto L_{q}$  such that $J= L_{m}$. If we define $J^{\prime} = L_{n}$, then $J^{\prime} \in \mathcal{B}(\mathcal{H})$ is anti self-adjoint and unitary such that $JJ^{\prime} = -J^{\prime}J$.
	\end{note}
	\begin{defn}\label{Definitioncf}
		Let $m,n \in \mathbb{S}$ with $mn = -nm$. Let $T \in\mathcal{B}(\mathcal{H})$ be normal and  $T_{+} \in \mathcal{B}(\mathcal{H}^{Jm}_{+})$ be such that $\widetilde{T_{+}} = T$.  If $f \in {C}(\sigma_{S}(T), \mathbb{H})$, then by Lemma \ref{decomposition}, there exist  unique $F_{1}, F_{2} \in {C}(\sigma_{S}(T), \mathbb{C}_{m})$ such that
		\begin{equation*}
		f(q) = F_{1}(q) + F_{2}(q) \cdot n,\text{for all } \; q \in \mathbb{H}.
		\end{equation*}
		Define
		\begin{equation}\label{defofcontF_calculus}
		f(T) = F_{1}(T) + F_{2}(T) J^{\prime},
		\end{equation}
		where $J^{\prime}$ is defined as in Note \ref{noteforJ}.
		Since $F_{1},F_{2} \in {C}(\sigma_{S}(T), \mathbb{C}_{m}) \subseteq C_{\mathbb{C}_{m}}((\sigma_{S}(T), \mathbb{H}))$, by Lemma \ref{decomposition}, the operators $F_{1}(T)$ and $F_{2}(T)$ are well-defined. Thus $f(T)$ is well-defined.
	\end{defn}
	\begin{lemma}\label{f(T)commutesJprime}
		Let $T \in \mathcal{B}(\mathcal{H})$ be normal and $m \in \mathbb{S}$. If $f \in C(\sigma_{S}(T), \mathbb{H})$, then $f(T)J^{\prime} = J^{\prime}f(T)^{*}$, where $J^{\prime}$ is same as in Note \ref{noteforJ}.
	\end{lemma}
	\begin{proof}
		Define $g \colon \sigma_{S}(T) \to \mathbb{H}$ by
		\begin{equation*}
		g(q) = f(q) \cdot n, \; \text{for all}\; q \in \sigma_{S}(T).
		\end{equation*}
		By Definition \ref{Definitioncf}, we have $g(T) = f(T) J^{\prime}$.
		Since $mn = -nm$, we write $g(q) = n \cdot \overline{f(q)}$. This implies that
		\begin{equation*}
		g(T) = J^{\prime} f(T)^{*}.
		\end{equation*}
		Thus $f(T)J^{\prime} = J^{\prime}f(T)^{*}$.
	\end{proof}
	It is worth to note that our method works for non slice continuous functions also. Here we illustrate this fact.
	
	We recall the definition of slice function. Let $m \in \mathbb{S}$ and $\mathcal{K}$ be a non-empty subset of $\mathbb{C}_{m}$. The complexification of $\mathbb{H}$ is given by
	\begin{equation*}
	\mathbb{H}_{\mathbb{C}_{m}}:= \mathbb{H} \otimes_{\mathbb{R}}\mathbb{C}_{m}.
	\end{equation*}
	If $w \in \mathbb{H}_{\mathbb{C}_{m}}$, then $w = q + m p$, for some $q,p \in \mathbb{H}$.
	\begin{defn}\cite[Defination 6.1]{Ghiloni}
		Let ${S} \ \colon \mathcal{K} \to \mathbb{H}_{\mathbb{C}_{m}}$ be a map and  ${S}_{1}, {S}_{2} \colon \mathcal{K} \to \mathbb{H}$ be components of ${S}$ that is ${S}(z) = {S}_{1}(z)+ m\; {S}_{2}(z)$, for all $z \in \mathcal{K}$. We say that ${S}$ is a stem function on $\mathcal{K}$ if ${S}_{1}(\overline{z}) = {S}_{1}(z)$ and ${S}_{2}(\overline{z}) = - {S}_{2}(z)$, for every $z \in \mathcal{K}$.
		
		Note that if ${S}_{1}$ and $S_{2}$ are continuous, then ${S}$ is continuous.
	\end{defn}
	\begin{defn} \cite[Definition 6.3]{Ghiloni}
		Every stem function  ${S} = {S}_{1} + i {S}_{2}$  induces a slice function $\mathcal{I}({S}) \colon \Omega_{\mathcal{K}} \longrightarrow \mathbb{H} $ on $\Omega_{\mathcal{K}}$ as follows:
		
		if $q = \alpha + j \beta \in \Omega_{\mathcal{K}}$ for some $\alpha, \beta \in \mathbb{R}$ and $j \in \mathbb{S}$, then
		\begin{equation*}
		\mathcal{I}({S})(q):= {S}_{1}(z) + j {S}_{2}(z), \; \text{where } \; z = \alpha + i \beta \in \mathcal{K}.
		\end{equation*}
		If ${S}$ is continuous, then $\mathcal{I}({S})$ is a continuous slice function on $\Omega_{\mathcal{K}}$.
	\end{defn}
	
	\begin{eg}\label{eg1}
		Let $\mu$ be a positive $\sigma$ - additive  measure on $\mathbb{S}$. Let  $m,n \in \mathbb{S}$ with $mn = -nm$.  Then $T \colon L^{2}(\mathbb{S}; \mathbb{H};\mu) \to L^{2}\big(\mathbb{S}; \mathbb{H};\mu \big)$ defined by
		\begin{equation*}
		T(g)(s) = sg(s), \; \text{for all}\; g \in L^{2}\big(\mathbb{S}; \mathbb{H};\mu\big), s \in \mathbb{S},
		\end{equation*}
		is a normal operator. The adjoint of $T$ is given by
		\begin{equation*}
		T^{*}(g)(s)= \overline{s}g(s) = -s g(s) = -T(g)(s).
		\end{equation*}
		Moreover, $T^{*}T(g)(s) = |s|^{2}g(s) = g(s) = TT^{*}(g)(s)$. This implies that $T$ is anti self-adjoint and unitary. Therefore,  $\sigma_{S}(T) = \mathbb{S}$.
		
		Let us take $J = T$. Then
		\begin{equation*}
		L^{2}(\mathbb{S}; \mathbb{H}; \mu)^{Jm}_{\pm}:= \big\{g \in L^{2}(\mathbb{S}; \mathbb{H}; \mu) : sg(s) = \pm g(s)m, \; \text{for all}\; s \in \mathbb{S}\; \big\}.
		\end{equation*}
		Then  $T_{+} \colon L^{2}(\mathbb{S}; \mathbb{H};\mu)^{Jm}_{+} \to L^{2}(\mathbb{S}; \mathbb{H};\mu)^{Jm}_{+}$  defined by
		\begin{equation*}
		T_{+}(g)(s) = T(g)(s)=s g(s), \; \text{for all}\; g \in L^{2}(\mathbb{S}; \mathbb{H};\mu)^{Jm}_{+}
		\end{equation*}
		is a $\mathbb{C}_{m}$- linear operator such that $\widetilde{T_{+}}=T$.  By Remark \ref{rmkforspT+}, we have
		\begin{equation*}
		\sigma(T_{+}) = \sigma_{S}(T)\cap \mathbb{C}_{m}^{+} = \mathbb{S}\cap \mathbb{C}_{m}^{+} = \{m\}.
		\end{equation*}
		
		Now we discuss continuous functional calculus for $T$.	Let $f \colon \mathbb{S} \to \mathbb{H}$ defined by
		\begin{equation*}
		f(\alpha + j \beta) = (\alpha + m \beta)  + j {(\alpha -  m \beta)}, \; \text{for all}\; \alpha + j \beta \in \mathbb{S}.
		\end{equation*}
		We show that $f$ is not a slice function. Assume that  $f$ is a slice function, then there exist a stem function ${S}$ on $\sigma(T_{+})$ such that
		\begin{equation*}
		f(\alpha + j \beta) = \mathcal{I}({S})(\alpha + j \beta) = {S}_{1}(\alpha + m \beta) + j {S}_{2}(\alpha + m \beta)
		\end{equation*}
		This implies that ${S}_{1}(\alpha + m \beta) = (\alpha + m \beta) $ and ${S}_{2}(\alpha + m \beta) = (\alpha- m\beta)$. But
		\begin{equation*}
		{S}_{1}(\overline{\alpha + m \beta}) = (\alpha - m \beta)  \neq {S}_{1}(\alpha + m \beta).
		\end{equation*}
		It shows that ${S}$ is not a stem function, a contradiction. Thus $f$ is not a slice function.
		
		We show that $f$ is continuous. Let $\big(\alpha_{\ell} +j_{\ell} \beta_{\ell}\big) \subseteq \sigma_{S}(T) $ converges to $\alpha + j \beta$ as $\ell \to \infty$. That is
		\begin{equation*}
		\|(\alpha_{\ell} - \alpha)+ (j_{\ell}\beta_{\ell}-j\beta)\|^{2} \longrightarrow 0, \; \text{as }\; \ell \to \infty.
		\end{equation*}
		This implies that $\alpha_{\ell} \to \alpha$ and $j_{\ell}\beta_{\ell} \to j\beta$, as $\ell \to \infty$. In fact,
		\begin{equation*}
		\beta_{\ell} = |\beta_{\ell}| = |j_{\ell}\beta_{\ell}| \longrightarrow |j\beta| = |\beta| = \beta, \; \text{as} \; \ell \to \infty.
		\end{equation*}
		Then
		$f(\alpha_{\ell} + j_{\ell} \beta_{\ell}) = (\alpha_{\ell} + m \beta_{\ell}) + j (\alpha_{\ell}-m \beta) $ converges to $(\alpha + m \beta) + j (\alpha - m\beta)$, as $\ell \to \infty$. This shows that $f(\alpha_{\ell} + j_{\ell} \beta_{\ell})$ converges to  $f(\alpha + j \beta)$, as $\ell \to \infty$, and hence $f \in C(\sigma_{S}(T), \mathbb{H})$.
		
		Now we decompose $f$ as in Lemma \ref{decomposition}. Since $j \in \mathbb{S}$, we have $j = j_{1}m + j_{2}n + j_{3}mn$, for some $j_{\ell} \in \mathbb{R}, \ell = 1,2,3$. Define  $F_{1}$ and $F_{2}$ on $\sigma_{S}(T)$ by
		\begin{align*}
		F_{1}(\alpha + j \beta) &= (\alpha + m \beta) + j_{1}m (\alpha - m \beta) , \\
		F_{2}(\alpha + j \beta) &=  (j_{2}+j_{3}m)(\alpha + m \beta).
		\end{align*}
		Clearly, $F_{1}, F_{2} \in {C}(\mathbb{S}, \mathbb{C}_{m})\subseteq {C}_{\mathbb{C}_{m}}(\mathbb{S}, \mathbb{H})$ and
		\begin{equation*}
		f(\alpha + j \beta) = F_{1}(\alpha + j \beta) + F_{2}(\alpha + j \beta)\cdot n, \; \text{for all}\;\; \alpha + j \beta \in \mathbb{S}.
		\end{equation*}
		If ${F_{1}}_{+}:= F_{1}|_{\sigma(T_{+})}$ and $ {F_{2}}_{+}\colon = F_{2}|_{\sigma(T_{+})} $, then
		\begin{align*}
		{F_{1}}_{+}(m) = m +1  \;\text{and}\; {F_{2}}_{+}( m)= 0.
		\end{align*}
		Since ${F_{1}}_{+},\; {F_{2}}_{+} \in C(\sigma(T_{+}), \mathbb{C}_{m})$, then by continuous functional calculus in complex Hilbert spaces, we have
		\begin{equation*}
		{F_{1}}_{+}(T_{+}) = T_{+} + I_{+} \; \text{and}\; {F_{2}}_{+}(T_{+})= 0,
		\end{equation*}
		where $I_{+}$ is the identity operator on $H^{Jm}_{+}$.
		
		By Definition \ref{Definitioncmf}, we have
		\begin{align*}
		F_{1}(T)(g)(s)&= T(g)(s)+I(g)(s)\\
		&= s g(s) + g(s)\\
		&= (s+1)g(s).\\
		F_{2}(T)(g)(s)&=  0, \; \text{for all}\; g \in L^{2}(\mathbb{S};\mathbb{H};\mu).
		\end{align*}
		By  Definition \ref{Definitioncf} and Lemma \ref{f(T)commutesJprime},  the operator $f(T)$ is defined as follows:
		\begin{align*}
		f(T)(g)(s) &=  F_{1}(T)(g)(s) + F_{2}(T) J^{\prime}(g)(s) \\
		&= (s+1)g(s),
		\end{align*}
		for all $g \in L^{2}(\mathbb{S};\mathbb{H};\mu)$.
	\end{eg}
	\begin{rmk}
		In Example \ref{eg1}, $f(T)(g)(s) = (s+1)g(s)$, \text{for all} $g \in L^{2}(\mathbb{S};\mathbb{H};\mu)$. In particular,
		\begin{equation*}
		f(T)(g)(m) = f(m)g(m), \text{for all} \; g \in L^{2}(\mathbb{S};\mathbb{H};\mu).
		\end{equation*}
	\end{rmk}
	\begin{eg} Let $T$ be defined as in Example \ref{eg1}. Define $f \colon \mathbb{S} \to \mathbb{H}$ by
		\begin{equation*}
		f(\alpha + j \beta)  = e^{(\alpha+j \beta)} , \; \text{for all}\; \alpha+j \beta \in \mathbb{S}.
		\end{equation*}
		Let $m,n \in \mathbb{S}$ with $mn = -nm$. Then  $f \in C_{\mathbb{C}_{m}}(\mathbb{S}, \mathbb{H})$ is bounded. If $f_{+}:= f|_{\sigma(T_{+})}$, then $f_{+}(m) = e^{m}$. Define
		\begin{equation*}
		{p_{\ell}}_{+}( m) = \sum\limits_{t=0}^{\ell}\frac{ m^{t}}{t!}.
		\end{equation*}
		Then $f_{+}(m) = \lim\limits_{\ell \to \infty}{p_{\ell}}_{+}(m)$.
		
		By continuous functional calculus in complex Hilbert spaces, we have
		\begin{equation*}
		{p_{\ell}}_{+}(T_{+}) = \sum\limits_{t=0}^{\ell}\frac{T_{+}^{t}}{t!}.
		\end{equation*}
		By Remark \ref{boundedapprx}, we have
		\begin{align*}
		f(T)(g)(s) = \widetilde{f_{+}(T_{+})}(g_{+}+g_{-})(s) &= \lim\limits_{\ell \to \infty}\widetilde{{p_{\ell}}_{+}(T_{+})}(g_{+}+g_{-})(s)\\
		&= \lim\limits_{\ell \to \infty}\sum\limits_{t = 0}^{\ell}\frac{T^{t}}{t!}(g)(s)\\
		&=\lim\limits_{\ell \to \infty}\sum\limits_{t = 0}^{\ell}\frac{s^{t}}{t!}g(s)\\
		&= e^{s}g(s)\\
		&=f(s)g(s),
		\end{align*}
		for all $g \in L^{2}(\mathbb{S}; \mathbb{H}; \mu)$.
	\end{eg}
	\section{Integral Representation}

	In this section we establish integral representation of  bounded quaternionic normal operators. First we define quaternionic projection valued measure based on \cite[Definition 12.17]{Rudin}.
	
	\begin{defn}\label{qspmeas}
		Let $T\in \mathcal{B}(\mathcal{H})$ be normal and $T_{+} \in \mathcal{B}(\mathcal{H}^{Jm}_{+})$ such that $\widetilde{T_{+}}=T$. Let  $ \Sigma_{\sigma(T_{+})} $ be the $\sigma$- algebra of  $\sigma(T_{+})$. A quaternionic projection valued spectral measure on $\sigma(T_{+})$ is a map  $ F \colon \Sigma_{\sigma(T_{+})} \to \mathcal{B}(\mathcal{H}) $ satisfying the following properties:
		\begin{enumerate}
			\item $ F(\emptyset) = 0 $ and $ F({\sigma(T_{+})}) = I$.
			\item $F(\omega)^{*}=F(\omega)$, for all $ \omega \in \Sigma_{\sigma(T_{+})}$.
			\item If $ \omega_{1}, \omega_{2} \in \Sigma_{\sigma(T_{+})}, $ then $ F(\omega_{1}\cap \omega_{2}) = F(\omega_{1})\cdot F(\omega_{2})$.
			\item For $ x,y \in \mathcal{H},$ the map  $ F_{x,y} \colon \Sigma_{\sigma(T_{+})} \to \mathbb{H}$ defined by
			\begin{equation*}
			F_{x,y}(\omega) = \big\langle x \;\big|\; F(\omega)y\big\rangle, \; \text{for all } \; \omega \in \Sigma_{\sigma(T_{+})} ,
			\end{equation*}
			is a $\mathbb{H}$- valued measure on $\sigma(T_{+})$.
		\end{enumerate}
	\end{defn}
	\begin{rmk}
		If we consider $\sigma$- algebra of Borel subsets of $\sigma(T_{+})$, then it is customary to add the condition that $F_{x,y}$ should be regular Borel measure for each $x,y \in \mathcal{H}$.
	\end{rmk}
	\begin{note}
		The  measure defined in Definition \ref{qspmeas} is also called qPVM (See \cite[Section 1.2]{Ghiloni2} for details).
	\end{note}
	The meaning of the integral of  general bounded measurable $\mathbb{H}$-valued functions with respect to quaternionic projection valued measure is similar to \cite[Definition 3.11]{Ghiloni2}.
	
	%
	If $f \colon \sigma(T_{+}) \to \mathbb{H}$ is a simple function, then $f = \sum\limits_{\ell = 1}^{n} q_{\ell} \chi_{\omega_{\ell}}$, for some $q_{\ell} \in \mathbb{H}$ and $\omega_{\ell} \in \Sigma_{\sigma(T_{+})}$. Moreover,
	\begin{equation*}
	\int\limits_{\sigma(T_{+})} f dF = \sum\limits_{\ell = 1}^{n} q_{\ell} F(\omega_{\ell}).
	\end{equation*}
	%
	\begin{defn} \cite[Definition 3.11]{Ghiloni2}
		Given $\omega \in \Sigma_{\sigma(T_{+})}$ and a bounded measurable function $f \colon \sigma(T_{+}) \to \mathbb{H}$, the integral of $f$ with respect to $F$ is the operator in $\mathcal{B}(H)$ defined by the following limit:
		\begin{equation*}
		\int\limits_{\omega} f dF := \lim\limits_{\ell \to \infty} \int\limits_{\sigma(T_{+})} f_{\ell} dF,
		\end{equation*}
		where $\{f_{\ell}\}_{\ell \in \mathbb{N}}$ is any sequence of simple functions on $\sigma(T_{+})$ such that $\|f_{\ell}-  \widetilde{f}\|_{\infty} \to 0$, as $\ell \to \infty$ and $\widetilde{f}$ extends $f$ to the null function outside $\omega$.
	\end{defn}
	
	Now we prove the existence of quaternionic projection valued measure.
	\begin{prop}\label{existenceofqpvm}
		Let $T \in \mathcal{B}(\mathcal{H})$ be normal and $m,n \in \mathbb{S}$ with $mn = -nm$. Let $T_{+} \in \mathcal{B}(\mathcal{H}^{Jm}_{+})$ be such that $\widetilde{T_{+}}=T$.  If $E$ is the complex projection valued spectral measure on $\sigma(T_{+})$, then the mapping $F \colon \sum_{\sigma({T_{+}})} \to  \mathcal{B}(\mathcal{H})$ defined by
		\begin{equation}\label{eq:qpvm}
		F(\omega) = \widetilde{E(\omega)},
		\end{equation}
		is a quaternion projection valued measure.
	\end{prop}
	\begin{proof} We show that $F$ satisfies all the properties listed in Definition \ref{qspmeas}. Since $E(\omega) \in \mathcal{B}(\mathcal{H}^{Jm}_{+})$ is a projection and by Equation (\ref{eq:qpvm}) and Proposition \ref{extension1}, it is clear that $F(\omega) \in \mathcal{B}(\mathcal{H})$ is a projection, and $JF(\omega) = F(\omega)J$, for all $\omega \in \sum_{\sigma(T_{+})}$.	If $\omega_{1}, \omega_{2} \in \sum_{\sigma(T_{+})}$, then
		\begin{equation*}
		F(\omega_{1} \cap \omega_{2}) = \widetilde{E(\omega_{1} \cap \omega_{2})} = \widetilde{E(\omega_{1})} \widetilde{E(\omega_{2})} = F(\omega_{1}) F(\omega_{2}).
		\end{equation*}
		Let $x,y \in \mathcal{H}$. Then $x = x_{1} + x_{2} \cdot n$ and $y = y_{1} + y_{2} \cdot n$, for some $x_{\ell},y_{\ell} \in \mathcal{H}^{Jm}_{+},\; \ell = 1,2$. If $\omega \in \sum_{\sigma(T_{+})}$, then	
		\begin{align*}
		F_{x,y}(\omega) &= \left\langle x | F(\omega)y\right\rangle \\
		&= \left\langle x_{1} + x_{2} \cdot n | E(\omega)y_{1} + E(\omega)y_{2} \cdot n\right\rangle \\
		&= \left\langle x_{1} | E(\omega) y_{1}\right\rangle + \left\langle x_{1} | E(\omega)y_{2} \cdot n\right\rangle + \left\langle x_{2} \cdot n | E(\omega)y_{1}\right\rangle + \left\langle x_{2} \cdot n | E(\omega) y_{2} \cdot n\right\rangle \\
		&= [E_{x_{1}, y_{1}}(\omega) + E_{y_{2}, x_{2}}(\omega)] + [E_{x_{1}, y_{2}}(\omega) - E_{y_{1},x_{2}}(\omega)] \cdot n.
		\end{align*}
		This implies that  $F_{x,y}$ is $\mathbb{H}$- valued measure on $\sigma(T_{+})$,  for all $x,y \in \mathcal{H}$.
	\end{proof}
	Since $T_{+}$, defined as in Note \ref{noteforT+}, is a $\mathbb{C}_{m}$- linear bounded normal operator on $\mathcal{H}^{Jm}_{+}$, there exist unique spectral measure $E$ on $\sigma(T_{+})$ by \cite[Theorems 12.22]{Rudin}.  We recall the following result.
	\begin{thm}\cite[subsection 12.24]{Rudin}\label{thm:classical}
		If $E$ is the  spectral measure on $\sigma(T_{+})$ and $f_{+}\in C(\sigma(T_{+}), \mathbb{C}_{m})$, then
		\begin{equation*}
		\big\langle a \; \big|\; f_{+}(T_{+})b\big\rangle = \int\limits_{\sigma(T_{+})} f_{+}(\lambda) \; dE_{a,b}(\lambda), \; \text{for all}\; a,b \in \mathcal{H}^{Jm}_{+}.
		\end{equation*}
		Moreover, every projection $E(\omega)$ commutes with every $B \in \mathcal{B}(\mathcal{H}^{Jm}_{+})$ which commutes with $T_{+}$.
	\end{thm}
	\begin{lemma}\label{lemmaforspthm}
		Let $T \in \mathcal{B}(\mathcal{H})$ be normal and $m,n \in \mathbb{S}$ with $mn = -nm$. Let $T_{+} \in \mathcal{B}(\mathcal{H}^{Jm}_{+})$ be such that $\widetilde{T_{+}} = T$ and $E$ be the  spectral measure on $\sigma(T_{+})$ as in Theorem \ref{thm:classical}.  If $f \in \mathcal{C}_{\mathbb{C}_{m}}(\sigma_{S}(T), \mathbb{H})$, then
		\begin{equation*}
		\overline{n} \cdot \int\limits_{\sigma(T_{+})} f_{+} \ dE_{a,b} = - \int\limits_{\sigma(T_{+})} \overline{f_{+}} \ dE_{b,a} \ \cdot n, \; \text{for all} \; a,b \in\mathcal{H}^{Jm}_{+},
		\end{equation*}
		where $f_{+}$ is the restriction of $f$ onto $\sigma(T_{+})$.
	\end{lemma}
	\begin{proof} Let $f_{+}$ be same as in Lemma \ref{invconti} and $K=\sigma(T_{+})$. By the definition of integral in Theorem \ref{thm:classical}, we have
		\begin{align*}
		\overline{n} \cdot \int\limits_{\sigma(T_{+})} f_{+} \ dE_{a,b} = -{n} \cdot \big\langle a \;\big|\; f_{+}(T_{+})b\big\rangle &= - \big\langle f_{+}(T_{+})b \;\big|\; a\big\rangle \cdot {n} \\
		&= - \left\langle b \;\big| f_{+}(T_{+})^{*}a\right\rangle \cdot n \\
		&= - \int\limits_{\sigma(T_{+})} \overline{f_{+}} \ dE_{b,a} \ \cdot n.\qedhere
		\end{align*}
	\end{proof}
	\begin{thm}\label{quaternionicspthm}
		Let $T \in \mathcal{B}(\mathcal{H})$ be normal and $m,n \in \mathbb{S}$ such that $mn = -nm$. Let $T_{+} \in \mathcal{B}(\mathcal{H}^{Jm}_{+})$ be such that $\widetilde{T_{+}} = T$, then there exists a Hilbert basis $\mathcal{N}_{m}$  of $\mathcal{H}$, and a unique quaternion projection valued spectral  measure $F$ on Borel subsets of $\sigma(T_{+})$  such that
		\begin{enumerate}
			\item if  $f \in C_{\mathbb{C}_{m}}(\sigma_{S}(T), \mathbb{H})$, then
			
			\begin{equation}\label{quatspthmeq}
			\left\langle x \; \big|\; f(T)y\right\rangle = \int\limits_{\sigma(T_{+})} f(\lambda)\ dF_{x,y}(\lambda), \;\text{for all}\; x ,y \in \mathcal{H}.
			\end{equation}
			
			\item If $f \in {C}(\sigma_{S}(T), \mathbb{H})$, then there exist unique $F_{1},F_{2} \in {C}(\sigma_{S}(T), \mathbb{C}_{m})$ such that $f(q) = F_{1}(q) + F_{2}(q) \cdot n$, for all $q \in \sigma_{S}(T)$. Moreover,
			\begin{equation*}
			\left\langle x \;\big|\; f(T)y\right\rangle = \int\limits_{\sigma(T_{+})} F_{1}(\lambda) \ dF_{x,y}(\lambda) +  \int\limits_{\sigma(T_{+})} F_{2}(\lambda)  \ dF_{x, n \cdot y}(\lambda) ,
			\end{equation*}
			for all $x,y \in \mathcal{H}$.
			
			Here the left multiplication by `n' is induced by the Hilbert basis $\mathcal{N}_{m}$.
			\item Let  $S \in \mathcal{B}(\mathcal{H})$ be such that $ST = TS$ and $ST^{*} =T^{*}S$, then
			\begin{enumerate}
				\item $Sf(T) = f(T)S$, for every $f \in {C}_{\mathbb{C}_{m}}(\sigma_{S}(T), \mathbb{H})$.
				\item  $SF(\omega) = F(\omega)S$, for every Borel set $\omega$ of $\sigma(T_{+})$.
			\end{enumerate}
		\end{enumerate}
	\end{thm}
	\begin{proof}
		Proof of {\it{(1)}:} If $f \in C_{\mathbb{C}_{m}}(\sigma_{S}(T), \mathbb{H})$, then by Equation(\ref{defofF_calculus}), the operator $f(T)$ is defined as the unique extension of $f_{+}(T_{+})$.
		Since $T_{+}$ is a bounded $\mathbb{C}_{m}$- linear normal operator on $\mathcal{H}^{Jm}_{+}$, by Theorem \ref{thm:classical}, there exist  unique spectral measure $E$ on $\sigma(T_{+})$ such that
		\begin{equation*}
		\big\langle a \; \big|\; f_{+}(T_{+})b\big\rangle = \int\limits_{\sigma(T_{+})} f_{+}(\lambda) \; dE_{a,b}(\lambda), \; \text{for all}\; a,b \in \mathcal{H}^{Jm}_{+},
		\end{equation*}
		where $f_{+}=f|_{\sigma(T_{+})}$.
		Let $F$ be the quaternionic projection valued spectral measure on $\sigma(T_{+})$ as in Proposition \ref{existenceofqpvm}. Then for every  $\omega \in \sum_{\sigma(T_{+})}$, and $ x = x_{1}+x_{2}\cdot  n, \;\; y = y_{1}+ y_{2}\cdot n \in \mathcal{H}$,  we have
		\begin{align*}
		&\int\limits_{\omega} dF_{x,y}(\lambda) \\
		&=  \left\langle x_{1} + x_{2} \cdot n \;\big|\;E(\omega)y_{1} + E(\omega)y_{2} \cdot n \right\rangle\\
		&= \left\langle x_{1}\;\big|\;E(\omega)y_{1} \right\rangle + \left\langle x_{1} \;\big|\; E(\omega)y_{2} \right\rangle \cdot n + \overline{n} \ \left\langle x_{2} \;\big|\; E(\omega)y_{1} \right\rangle + \overline{n} \left\langle x_{2} \;\big|\; E(\omega)y_{2}\right\rangle \cdot n\\
		&= \int\limits_{\omega}dE_{x_{1},y_{1}}(\lambda) +\int\limits_{\omega}dE_{x_{1},y_{2}}(\lambda) \cdot n +\overline{n} \int\limits_{\omega}dE_{x_{2},y_{1}}(\lambda) +\overline{n} \int\limits_{\omega}dE_{x_{2},y_{2}}(\lambda) \cdot n.
		\end{align*}
		This implies that
		\begin{align*}
		\int\limits_{\sigma(T_{+})}f(\lambda) \ dF_{x,y}(\lambda)
		&= \int\limits_{\sigma(T_{+})} f_{+}(\lambda) \ dE_{x_{1},y_{1}}(\lambda) + \int\limits_{\sigma(T_{+})} f_{+}(\lambda) \ dE_{x_{1}, y_{2}}(\lambda) \cdot n \\
		& \hspace{0.8cm} + \overline{n} \int\limits_{\sigma(T_{+})}\  f_{+}(\lambda) \ dE_{x_{2}, y_{1}} (\lambda) + \overline{n} \int\limits_{\sigma(T_{+})} f_{+}(\lambda) \ dE_{x_{2},y_{2}}(\lambda) \cdot n.
		\end{align*}
		By Lemma \ref{lemmaforspthm}, we have
		\begin{align*}
		\int\limits_{\sigma(T_{+})}&f(\lambda) \ dF_{x,y}(\lambda)\\
		&= \int\limits_{\sigma(T_{+})} f_{+}(\lambda) \ dE_{x_{1},y_{1}}(\lambda) + \int\limits_{\sigma(T_{+})} f_{+}(\lambda) \ dE_{x_{1}, y_{2}}(\lambda) \cdot n \\
		& \hspace{0.8cm} - \int\limits_{\sigma(T_{+})}\overline{f_{+}}(\lambda) \ dE_{y_{1},x_{2}} \cdot n + \int\limits_{\sigma(T_{+})}\overline{f_{+}}(\lambda) \ dE_{y_{2},x_{2}} \\
		& = \left\langle x_{1} \;\big|\; f_{+}(T_{+})y_{1}\right\rangle + \left\langle x_{1}\;\big|\;f_{+}(T_{+})y_{2}\right\rangle \cdot n - \left\langle y_{1}\;\big|\;f_{+}(T_{+})^{*}x_{2}\right\rangle \cdot n  \\
		& \hspace{0.8cm } + \left\langle y_{2}\;\big|\; f_{+}(T_{+})^{*}x_{2} \right\rangle\\
		&= \left\langle x_{1} \;\big|\; f_{+}(T_{+})y_{1}\right\rangle + \left\langle x_{1} \;\big|\; f_{+}(T_{+})y_{2}\right\rangle \cdot n + \overline{n} \left\langle x_{2} \;\big|\; f_{+}(T_{+})y_{1}\right\rangle \\
		& \hspace{0.8cm } + \overline{n}\left\langle x_{2} \;\big|\; f_{+}(T_{+})y_{2}\right\rangle \cdot n \\
		&= \left\langle x_{1} \;\big|\; \widetilde{f_{+}(T_{+})}y\right\rangle + \left\langle x_{2} \cdot n \;\big|\; \widetilde{f_{+}(T_{+})}y \right\rangle\\
		&= \left\langle x \;\big|\; f(T)y\right\rangle.
		\end{align*}
		Therefore
		\begin{equation*}
		\left\langle x \;\big|\; f(T)y\right\rangle = \int\limits_{\sigma_{S}(T)\cap \mathbb{C}_{m}^{+}} f(\lambda) \ dF_{x,y}(\lambda).
		\end{equation*}
		
		Proof of {\it(2):} If $f \in {C}(\sigma_{S}(T), \mathbb{H})$, then by Lemma \ref{decomposition}, there exist $F_{1}, F_{2} \in {C}(\sigma_{S}(T), \mathbb{C}_{m})$ such that $f(q) = F_{1}(q) + F_{2}(q) \cdot n$, for all $q \in \sigma_{S}(T)$. By Definition (\ref{Definitioncf}), we have
		\begin{align*}
		\left\langle x | f(T)y\right\rangle &= \left\langle x \;\big|\; F_{1}(T)y  \right\rangle + \left\langle x \;\big|\; F_{2}(T)J^{\prime}y\right\rangle \\
		&= \int\limits_{\sigma(T_{+}) } F_{1}(\lambda) \ dF_{x,y}(\lambda) + \int\limits_{\sigma(T_{+}) } F_{2}(\lambda)\ dF_{x,J^{\prime}y}(\lambda) \\
		&= \int\limits_{\sigma(T_{+}) } F_{1}(\lambda) dF_{x,y}(\lambda) + \int\limits_{\sigma(T_{+}) } F_{2}(\lambda) \ dF_{x, n \cdot y}(\lambda).
		\end{align*}
		
		Here $J^{\prime} = L_{n}$, the left multiplication induced by the Hilbert basis $\mathcal{N}_{m}$.	
		Proof of {\it{(3)(a):}} Since $S$ commutes with both $T$ and $T^{*}$, it is clear that $S$ commutes with $T-T^{*}$. The construction of $J$ in \cite[Theorem 5.9]{Ghiloni} shows that $J$ commutes with $S$.  Then there exists  a bounded $\mathbb{C}_{m}$- linear operator $S_{+}$ on $\mathcal{H}^{Jm}_{+}$  such that $S = \widetilde{S_{+}}$, and $S_{+}T_{+} = T_{+}S_{+}$. By Theorem \ref{thm:classical}, we have  $S_{+}f_{+}(T_{+}) = f_{+}(T_{+})S_{+}$, for all $f \in C_{\mathbb{C}_{m}}(\sigma_{S}(T) , \mathbb{H})$. Thus
		\begin{equation*}
		Sf(T) = f(T)S,\; \text{for all}\; f \in C_{\mathbb{C}_{m}}(\sigma_{S}(T) , \mathbb{H}).
		\end{equation*}
		\noindent Proof of {\it{(3)(b):}} Since $S$ commutes with $J$ and by Theorem \ref{thm:classical}, we have  $S_{+}E(\omega) = E(\omega)S_{+}$, for every Borel subset $\omega$ of $\sigma(T_{+})$. This implies that
		\begin{equation}\label{bddcommutes}
		SF(\omega) = \widetilde{S_{+}} \widetilde{E(\omega)} = \widetilde{S_{+}E(\omega)} = \widetilde{E(\omega)S_{+}} = F(\omega) S. 	
		\end{equation}
		Now we show the uniqueness. Suppose $F$ and $G$  are quaternionic projection valued spectral measures satisfying Equation (\ref{quatspthmeq}). Since $J$ commutes with both $T$ and $T^{*}$, by Equation(\ref{bddcommutes}), we have
		\begin{equation*}
		F(\omega)J = JF(\omega) \; ; \; JG(\omega) = G(\omega)J, \; \text{for all} \; \omega \in \Sigma_{\sigma(T_{+})}.
		\end{equation*}
		By Proposition \ref{extension1}, for every $\omega \in \Sigma_{\sigma(T_{+})}$, there exist a unique $\mathbb{C}_{m}$- linear operators $F(\omega)_{+}$ and $G(\omega)_{+}$  on $H^{Jm}_{+}$ such that $	
		F(\omega) = \widetilde{F(\omega)_{+}} \; \text{and}\; G(\omega) = \widetilde{G(\omega)_{+}}$.
		
		Define $F_{+}$ and $G_{+}$ on $\Sigma_{\sigma(T_{+})}$ by
		\begin{equation*}
		F_{+}(\omega) = F(\omega)_{+} \; \text{and}\; G_{+}(\omega) = G(\omega)_{+}
		\end{equation*}
		respectively. By the Proposition \ref{extension1}, Both $F_{+}$ and $G_{+}$ are complex projection valued spectral measures on $\sigma(T_{+})$. Moreover, if $a,b \in \mathcal{H}^{Jm}_{+}$, then
		\begin{equation*}
		\left\langle a \;\big|\; f_{+}(T_{+})b \right\rangle = \int\limits_{\sigma(T_{+})} f_{+}(\lambda) \ dF_{+_{a,b}}(\lambda) = \int\limits_{\sigma(T_{+})}  f_{+}(\lambda) \ dG_{+_{a,b}}(\lambda).
		\end{equation*}
		By the uniqueness of complex projection valued spectral measure on $\sigma(T_{+})$ as in Theorem \ref{thm:classical}, we have  $F_{+} = G_{+}$. This shows that
		\begin{equation*}
		F(\omega) = \widetilde{F_{+}(\omega)} = \widetilde{G_{+}(\omega)} = G(\omega), \; \text{for all}\; \omega \in \Sigma_{\sigma(T_{+})}.\qedhere
		\end{equation*}
	\end{proof}
	\section*{Acknowledgment}
	The second author is thankful to INSPIRE (DST) for the support in the form of fellowship (No. DST/ INSPIRE Fellowship/2012/IF120551), Govt of India.

\end{document}